\documentclass[12pt,]{article}
\usepackage{amsmath}
\usepackage{amscd}

\usepackage{amssymb}
\usepackage{latexsym}
\newtheorem{theorem}{Theorem}[section]
\newtheorem{lemma}[theorem]{Lemma}
\newtheorem{corollary}[theorem]{Corollary}
\newtheorem{proposition}[theorem]{Proposition}

\newtheorem{remark}[theorem]{Remark}
\newtheorem{definition}[theorem]{Definition}

\newenvironment{proof}{{\em Proof.}\ \ \ }{\unskip\nobreak\hfil\penalty50
\hskip1em\hbox{}\nobreak\hfil $\Box$
\parfillskip=0pt \finalhyphendemerits=0 \par\medskip\noindent}
\newcommand{\bgl}{\begin{equation}}         
\newcommand{\egl}{\end{equation}}
\newcommand{\bgln}{\begin{eqnarray}}        
\newcommand{\egln}{\end{eqnarray}}
\newcommand{\bglnoz}{\begin{eqnarray*}}     
\newcommand{\eglnoz}{\end{eqnarray*}}
\newcommand{\btheo}{\begin{theorem}}
\newcommand{\etheo}{\end{theorem}}
\newcommand{\blemma}{\begin{lemma}}
\newcommand{\elemma}{\end{lemma}}
\newcommand{\bproof}{\begin{proof}}
\newcommand{\eproof}{\end{proof}}
\newcommand{\bbew}{\begin{beweis}}
\newcommand{\ebew}{\end{beweis}}
\newcommand{\bremark}{\begin{remark}\em}
\newcommand{\eremark}{\end{remark}}
\newcommand{\bdefin}{\begin{definition}}
\newcommand{\edefin}{\end{definition}}
\newcommand{\bprop}{\begin{proposition}}
\newcommand{\eprop}{\end{proposition}}
\newcommand{\bcor}{\begin{corollary}}
\newcommand{\ecor}{\end{corollary}}
\newcommand{\n}{\par\noindent}

\newcommand{\mn}{\par\medskip\noindent}

\def\SEMI{\mbox{$\times\kern-2pt\vrule height5pt width.6pt \kern3pt $}}

\def\Cz{\mathbb{C}}

\def\Nz{\mathbb{N}}

\def\Rz{\mathbb{R}}
\def\Zz{\mathbb{Z}}

\begin{document}
\title{Cyclic Theory and the Bivariant Chern-Connes Character}%
\author{Joachim Cuntz \\ Mathematisches Institut \\ Universit\"at
M\"unster
\\ Einsteinstr. 62 \\ 48149 M\"unster \\ cuntz@math.uni-muenster.de}%
\thanks{Research supported by the Deutsche Forschungsgemeinschaft}
\maketitle
\begin{abstract}\noindent
We give a survey of cyclic homology/cohomology theory including
a detailed discussion of cyclic theories for various classes of
topological algebras. We show how to associate cyclic classes
with Fredholm modules and $K$-theory classes and how to
construct a completely general bivariant Chern-Connes character
from bivariant $K$-theory to bivariant cyclic theory.
\end{abstract}
\section{Introduction}
The two fundamental ``machines" of non-commutative geometry are
cyclic homology and (bivariant) topological $K$-theory. In the
present notes we describe these two theories and their
connections. Cyclic theory can be viewed as a far reaching
generalization of the classical de Rham cohomology, while
bivariant $K$-theory includes the topological $K$-theory of
Atiyah-Hirzebruch as a very special case. \mn The classical
commutative theories can be extended to a striking amount of
generality. It is important to note however that the new
theories are by no means based simply on generalizations of the
existing classical constructions. In fact, the constructions are
quite different and give, in the commutative case, a new
approach and an unexpected interpretation of the well-known
classical theories. One aspect is that some of the properties
of the two theories become visible only in the non-commutative
category. For instance both theories have certain universality
properties in this setting. \mn Bivariant $K$-theory has first
been defined and developed by Kasparov on the category of
$C^*$-algebras (possibly with the action of a locally compact
group) thereby unifying and decisively extending previous work
by Atiyah-Hirzebruch, Brown-Douglas-Fillmore and others.
Kasparov also applied his bivariant theory to obtain striking
positive results on the Novikov conjecture. Very recently, it
was discovered that in fact, bivariant topological $K$-theories
can be defined on a wide variety of topological algebras ranging
from discrete algebras and very general locally convex algebras
to e.g. Banach algebras or $C^*$-algebras (possibly equipped
with a group action). If $E$ is the covariant functor from such
a category of algebras given by topological $K$-theory or also
by periodic cyclic homology, then it satisfies the following
three fundamental properties:
\begin{itemize}
\item [(E1)] $E$ is diffeotopy invariant, i.e., the evaluation map $ev_t$
in any point $t\in[0,1]$ induces an isomorphism $E(ev_t):E(\frak
A[0,1])\to E(\frak A)$ for any $\frak A$ in $C$. Here $\frak
A[0,1]$ denotes the algebra of $\frak A$-valued
$C^\infty$-functions on $[0,1]$.
\item [(E2)] $E$ is stable, i.e., the canonical inclusion
$\iota:\frak A\to\frak K\hat{\otimes}\frak A$, where $\frak K$
denotes the algebra of infinite $\Nz\times\Nz$-matrices with
rapidly decreasing coefficients, induces an isomorphism
$E(\iota)$ for any $\frak A$ in $C$.
\item [(E3)] $E$ is half-exact, i.e., each extension $0\to\frak I
\to\frak A\to\frak B\to 0$ in $C$ admitting a continuous linear
splitting induces a short exact sequence $E(\frak I) \to E(\frak
A)\to E(\frak B)$
\end{itemize}
``Diffeotopy", i.e., differentiable homotopy is used in (E1) for
technical reasons in connection with the homotopy invariance
properties of cyclic homology. For $K$-theory, it could also be
replaced by ordinary continuous homotopy. \mn It turns out that
the bivariant $K$-functor is the universal functor from the
given category $C$ of algebras into an additive category $D$
(i.e., the morphism sets $D(\frak A,\frak B)$ are abelian
groups) satisfying these three properties \ref{univfunct}. \mn
Cyclic theory is a homology theory that has been developed,
starting from $K$-theory, independently by Connes and Tsygan.
Connes' construction was in fact directly motivated by
Kasparov's formalism for bivariant $K$-theory and in particular
for $K$-homology. A crucial role is played by so-called
Fredholm modules or spectral triples. Also Tsygan's work is
closely related to $K$-theory, \cite{Tsyg}. In fact, in his
approach, the new theory was  originally called ``additive
$K$-theory" and he pointed out that it is an additive version
of Quillen's definition of algebraic $K$-theory. It was
immediately realized that cyclic homology has close connections
with de Rham theory, Lie algebra homology, group cohomology and
index theorems. \mn The theory with the really good properties
is the periodic theory introduced by Connes. Periodic cyclic
homology $HP_*$ satisfies the three properties (E1), (E2),
(E3), \cite{CoIHES}, \cite{Good}, \cite{Wo1}, \cite{Wo2},
\cite{CQInv}. Combining this fact with the universality
property of bivariant $K$-theory leads to a multiplicative
transformation (the bivariant Chern-Connes character) from
bivariant $K$-theory to bivariant periodic cyclic theory. This
transformation is a vast generalization of the classical Chern
character in differential geometry. \mn The principal aim of
this volume is to give an account of cyclic theory. Here,
everything works in parallel algebraically as well as for
locally convex algebras (to name just two examples think of the
algebra of $C^\infty$-functions on a smooth manifold or of
algebras of pseudodifferential operators). Cyclic theory can be
introduced using rather different complexes, each one of them
having its own special virtues. Specifically, we will use the
following complexes or bicomplexes
\begin{itemize}
\item The cyclic bicomplex with various realizations:
\begin{eqnarray*} CC^n(A)\quad \mbox{with boundary operators}
\quad b,b', Q, \mbox{1}-\lambda
\\
\overline{CC}^n(\widetilde{A})\quad \mbox{with boundary
operator}\quad B-b
\\
\Omega(A)\quad \mbox{with boundary operator}\quad B-b
\end{eqnarray*}
The cyclic bicomplex is well suited
for the periodic theory as well as for the $\Zz$-graded
ordinary theory and for the connections between both.
\item The Connes complex $C^n_\lambda$
\\ It has the advantage, that concrete \emph{finite-dimensional}
cocycles often arise naturally as elements of $C^n_\lambda$.
The connection with Hochschild cohomology also fits naturally
into this picture.
\item The $X$-complex of any complete quasi-free extension of
the given algebra $A$. This complex is very useful for a
conceptual explanation of the properties of periodic theory, in
particular for proving excision, for the connections with
topological $K$-theory and the bivariant Chern-Connes
character. It also is the natural framework for all
infinite-dimensional versions of cyclic homology (analytic and
entire as well as asymptotic and local theory).
\end{itemize}
In the following sections we discuss the basic properties of
cyclic theory. We note however, that it is quite difficult to be
exhaustive and we don't even try to give a complete account of
all aspects of cyclic theory. For instance, an important notion
which we don't treat is the one of a cyclic object. But also
other important aspects have to be omitted. We focus on those
notions and results which, we think, are most relevant for
non-commutative geometry, including homotopy invariance, Morita
invariance, excision but also explicit formulas for the Chern
character associated to idempotents, invertibles, Fredholm
modules etc.. \mn After this we turn to a description of
bivariant $K$-theory and to the construction of the bivariant
Chern-Connes character, which generalizes the Chern character
for idempotents, invertibles or Fredholm modules mentioned
before. \mn As we pointed out already above, cyclic homology and
bivariant $K$-theory can be defined on different categories of
algebras - purely algebraically for algebras (over $\Rz$ or
$\Cz$) or on categories of topological algebras like locally
convex algebras, Banach algebras or C*-algebras. There are
different variants of the two theories which are adapted to the
different categories. In this text we treat cyclic theory in
the purely algebraic case (however restricting to algebras over
a field of characteristic 0) on the one hand. Concerning cyclic
theory for topological algebras on the other hand, we have to
make a choice. For the classical (``finite-dimensional") cyclic
theories we concentrate on the category of what we call
$m$-algebras. These are particularly nice locally convex
algebras (projective limits of Banach algebras). Their advantage
is that, both, ordinary cylic homology and bivariant topological
$K$-theory make perfect sense and that the bivariant
Chern-Connes-character can be constructed nicely on this
category. It is important to note however that the restriction
to the category of $m$-algebras is for convenience mainly. \n
Cyclic theory as well as bivariant $K$-theory can also be
treated on many other categories of algebras. In particular,
for the categories of Banach algebras or $C^*$-algebras there
are special variants of cyclic theory, namely the entire and
the local cyclic theory, which are especially designed for these
categories. An interesting new feature here is the existence of
\emph{infinite-dimensional} cohomology classes. These theories,
as well as a bivariant character from Kasparov's $KK$-theory to
the local cyclic theory, are discussed in sections 22 and 23
which have been written by Ralf Meyer and Michael Puschnigg.\mn
The text starts with a collection of examples of algebras,
locally convex algebras and certain extensions of algebras.
These serve as a reference for later and may be omitted at a
first reading. \mn \textbf{A copy of the complete article (68
pages) can be obtained as a SFB-preprint from our homepage at
wwwmath.uni-muenster.de}
\end{document}